\documentclass{gtart_a}
\pdfoutput=1

\title{Virtually Haken fillings and semi-bundles}

\author{Daryl Cooper}
\givenname{Daryl}
\surname{Cooper}
\address{Math Department\\UCSB\\Santa Barbara, CA 93106\\USA}
\email{cooper@math.ucsb.edu}
\urladdr{}

\author{Genevieve S Walsh}
\givenname{Genevieve S}
\surname{Walsh}
\address{{\rm [GW]}\qua Department of Math\\Tufts University\\Medford, MA 02155
\\USA\\and\\\newline
D\'epartement de Math\'ematiques\\UQAM\\Montr\'eal, QC H3C 3J7\\Canada}
\email{genevieve.walsh@tufts.edu}
\urladdr{}

\volumenumber{10}
\issuenumber{}
\publicationyear{2006}
\papernumber{49}
\lognumber{0721}
\startpage{2237}
\endpage{2245}

\doi{}
\MR{}
\Zbl{}

\arxivreference{math.gt/0407328}

\keyword{virtually Haken conjecture}
\keyword{Dehn filling}
\keyword{semi-bundles}
\subject{primary}{msc2000}{57M10}
\subject{secondary}{msc2000}{57M25}

\received{23 September 2004}
\revised{8 March 2006}
\accepted{24 October 2006}
\published{29 November 2006}
\publishedonline{29 November 2006}
\proposed{David Gabai}
\seconded{Cameron Gordon, Joan Birman}
\corresponding{}
\editor{CPR}
\version{}

\AtBeginDocument{\let\bar\wbar\let\tilde\wwtilde\let\hat\what}
\def\co{\mskip 0.7mu\colon\thinspace}
\makeop{orb}\makeop{Im}\makeop{incl}

\makeatletter
\def\cnewtheorem#1[#2]#3{\newtheorem{#1}{#3}[section]
\expandafter\let\csname c@#1\endcsname\c@theorem}

\newtheorem{theorem}{Theorem} \newtheorem{lemma}[theorem]{Lemma}
\cnewtheorem{corollary}[theorem]{Corollary}
\cnewtheorem{proposition}[theorem]{Proposition}

\theoremstyle{definition}

\makeatother  

\let\demo\proof

\newcommand{\zt}{\mathbb{Z}_2}

\begin{document}

\begin{asciiabstract}
Suppose that M is a fibered three-manifold whose fiber is a surface of
positive genus with one boundary component. Assume that M is not a
semi-bundle.  We show that infinitely many fillings of M along dM are
virtually Haken.  It follows that infinitely many Dehn-surgeries of
any non-trivial knot in the three-sphere are virtually Haken.
\end{asciiabstract}

\begin{htmlabstract}
Suppose that M is a fibered three-manifold whose fiber is a surface
of positive genus with one boundary component. Assume that M is not
a semi-bundle.  We show that infinitely many fillings of M along
&part; M are virtually Haken.  It follows that infinitely many
Dehn-surgeries of any non-trivial knot in the three-sphere are
virtually Haken.
\end{htmlabstract}

\begin{abstract} 
Suppose that $M$ is a fibered three-manifold whose fiber is a surface
of positive genus with one boundary component. Assume that $M$ is not
a semi-bundle.  We show that infinitely many fillings of $M$ along
$\partial M$ are virtually Haken.  It follows that infinitely many
Dehn-surgeries of any non-trivial knot in the three-sphere are
virtually Haken.
\end{abstract}

\maketitle

\section{Introduction}

In this paper {\em manifold} will always mean a compact, connected, orientable, possibly bounded, three-manifold. A {\em bundle} means a manifold which fibers over the circle. A {\em semi-bundle} is a manifold which is the union of two twisted $I$--bundles (over connected surfaces) whose intersection is the corresponding $\partial I$--bundle.  An irreducible, $\partial$--irreducible manifold  that contains a properly embedded incompressible surface is called {\it Haken}. A manifold is {\it virtually Haken} if  has a finite cover that is  Haken. 

Waldhausen's {\em virtually Haken conjecture} is that every
irreducible closed manifold with infinite fundamental group is
virtually Haken.  It was shown by Cooper and Long
\cite{CooperLongvirt} that {\em most} Dehn-fillings of an atoroidal
Haken manifold with torus boundary are virtually Haken provided the
manifold is not a bundle.

 \begin{theorem}\label{maintheorem} Suppose that $M$ is a bundle with fiber a compact surface $F$ and that $F$ has exactly one boundary component.  Also suppose that $M$ is not a  semi-bundle and not $S^1\times D^2.$ Then infinitely many Dehn-fillings of $M$ along $\partial M$ are virtually
Haken.\end{theorem}
\begin{corollary}\label{knot} Let $k$ be a knot in a homology three-sphere $N.$ Suppose that $N-k$ is irreducible and that $k$ does not bound a disk in $N.$ Then infinitely many Dehn-surgeries along $k$ are virtually
Haken. \end{corollary}

The main idea is to construct a surface of {\em invariant slope} (see \fullref{sec3}) in a particular finite cover of $M.$  Such surfaces are studied in arbitrary covers using representation theory in a sequel 
\cite{CooperWalsh2}. While writing this paper we noticed that Thurston's theory of bundles extends to semi-bundles, and in particular there are manifolds which are semi-bundles in infinitely many ways. We discuss this in the next section.

\medskip
We thank the referee for several helpful comments. The first author was partially supported by NSF grant DMS-0405963.  

\section{Bundles and semi-bundles}
 Various authors have studied semi-bundles, in particular Hempel and
Jaco
\cite{HempelJaco} and Zulli \cite{Z1,Z2}.
Suppose a manifold has a regular cover which is a surface bundle. 
We wish to
know when a particular fibration in the cover corresponds to a bundle or
semi-bundle structure on the quotient.  The following has the same flavor as some results of Hass~\cite{Hassmin}.

\begin{theorem} \label{fibering} Let $M$ be a compact, connected, orientable, irreducible
three-manifold,  $p\co \tilde M \rightarrow M$  a finite regular cover, and $G$ the
group of covering automorphisms.  Suppose that $\phi \co \tilde M \rightarrow S^1$ is a
fibration of $\tilde M$ over the circle. Suppose that the cyclic 
subgroup $V$ of
$H^1(\tilde M; {\Bbb Z})$ generated by $[\phi]$ is invariant under 
the action of
$G.$  Then one of the following occurs: 

\begin{enumerate} \item  The 
action of $G$ on $V$ is trivial. Then $M$ also fibers over the
circle. Moreover there is a fibering of $M$ which is covered by a
fibering of $\tilde{M}$ that is isotopic to the original fibering.
    
\item The action of $G$ on $V$ is non-trivial.  Then $M$ is a semi-bundle.
Moreover there is a semi-fibering of $M$ which is covered by a fibering of
$\tilde{M}$ that is isotopic to the original fibering.
\end{enumerate} \end{theorem}

\demo Define $N = \ker[\phi_*\co\pi_1\tilde{M}\rightarrow\pi_1S^1].$ Since $\phi$ is a fibration $N$ is finitely generated. If $N$ is cyclic then the fiber is a disc or annulus. In these cases the result is easy. Thus we may assume $N$ is not cyclic.
Because $V$ is $G$--invariant, it follows that $N$  is a normal subgroup of $\pi_1M$ and $Q = \pi_1M/N$ is infinite. 
Using \cite[Theorem 3]{HempelJaco} it follows that $M$ is a bundle or semi-bundle (depending on case 1 or 2) with fiber a compact surface $F$ and $N$ has finite index in $\pi_1F.$ The pull-back of this (semi)fibration of $M$ gives a fibration of $\tilde{M}$ in the cohomology class of $\phi$ and is therefore isotopic to the given fibration. \endproof

 Suppose that $G \cong (\zt)^n$ acts on a real vector space $V$ and let $X = Hom(G,{\Bbb C})$ denote the set of characters on $G.$ Then $X \cong Hom(G,  \zt).$ For each $\epsilon \in X$ there is a $G$--invariant {\em generalized $\epsilon$--eigenspace} $$V_{\epsilon} = \{\ v\in V\ :\ \ \forall g\in G\  \ g\cdot v = \epsilon(g)v\ \}.$$ Then $V$ is the direct sum of these subspaces $V_{\epsilon}.$  
 
Suppose that $M$ is an atoroidal irreducible manifold with boundary consisting of incompressible tori. According to Thurston there is a finite collection (possibly empty),  ${\cal C} = \{ C_1,\cdots,C_k\},$ called {\em fibered faces}. Each fibered face is  the interior of a certain top-dimensional face of the unit ball of the Thurston norm on $H_2(M,\partial M;{\Bbb R}).$ It is an
 open convex set  with the property that  fibrations of $M$ correspond to rational points in the projectivized space ${\Bbb P}(\cup_i C_i) \subset {\Bbb P}\left(H_2(M,\partial M;{\Bbb R})\right).$

Let $G= H_1(M;{\Bbb Z}/2).$ The regular cover $\tilde{M}_s$ of $M$ with covering group $G$ is called the {\em ${\Bbb Z}_2$--universal cover.}  Let ${\cal D} = \{ D_1,\cdots, D_l\}$ be the fibered faces for this cover. For each $\epsilon \in H^1(M;{\Bbb Z}_2)$ there is an $\epsilon$--eigenspace $H_{2,\epsilon}$ of $H_2({\tilde M}_s,\partial\tilde{M}_s;{\Bbb R}).$ For each 
$1\le i\le l$ and $\epsilon\in  H^1(M;{\Bbb Z}_2)$ we call $S_{i,\epsilon} = D_i\cap H_{2,\epsilon}$ a
{\em  semi-fibered face} if it is not empty. It is the interior of a compact convex polyhedron whose interior is in the interior of some fibered face for $\tilde{M}_s.$ Let $S_i$  be the union of the $S_{i, \epsilon}$ where $\epsilon$ is non-trivial. 
  
 \begin{theorem} \label{semifiberedfaces} With the above notation there is a bijection between isotopy classes of semi-fiberings of $M$ and rational points in ${\Bbb P}(\cup_i S_i).$ \end{theorem}

\demo A semi-fibration of $M$ gives such a rational point by considering the induced fibration on $\tilde{M}_s.$ The converse follows from \fullref{fibering}. We leave it as an exercise to check uniqueness up to isotopy.\endproof

\noindent We believe that  all points in ${\Bbb P}(\cup_i S_i)$ correspond to isotopy classes of non-transversally-orientable, transversally-measured, product-covered $2$--dimensional foliations of $M.$ This is true for rational points and therefore holds on a dense open set (using the fact that the set of non-degenerate twisted  $1$--forms is open). However, since we have no use for this fact, we have not tried very hard to prove it.

\medskip {\bf Definition}\qua A manifold is a {\em sesqui-bundle} if
it is both a
bundle and a semi-bundle.

\medskip An example is the torus bundle $M$ with monodromy
$-\mathrm{Id}$. This is the quotient of Euclidean three-space by the
group $\mathcal{G}_2$ (Wolf \cite[Theorem 3.5.5]{Wolf}).  $M$ has
infinitely many semi-fibrations with generic fiber a torus and two
Klein-bottle fibers. In addition, $M$ is a bundle thus a
sesqui-bundle.

A hyperbolic example may be obtained from $M$ as follows. Let $C$ be a 1--submanifold
in $M$ which is a small $C^1$--perturbation of a finite set of disjoint, immersed,
closed geodesics in
$M$ chosen so that:

\begin{enumerate}
\item No two components of $C$ cobound an annulus and no component
bounds a Mobius strip.

\item $C$ intersects every flat torus and flat Klein bottle.

\item Each component of $C$ is transverse to both a chosen fibration
  and semi-fibration.
\end{enumerate}

\noindent  Let $N$ be $M$ with a regular neighborhood of $C$ removed.  Then the
interior of $N$ admits a complete hyperbolic  metric. By (3)  it is a sesqui-bundle.
This answers a question of Zulli who asked in \cite{Z2} if  there are non-Seifert 
3--manifolds which are sesqui-bundles.

\section{Virtually Haken fillings}\label{sec3}

The following is well-known, but we include it here for ease of reference.
   \begin{lemma}\label{SFSfill} Suppose $M$ is Seifert fibered and has one boundary component . Then one of the following holds:

\begin{enumerate}
\item  $M$ is $D^2\times S^1$ or a twisted $I$--bundle over the Klein bottle. 
\item Infinitely many Dehn-fillings are virtually Haken. 
\end{enumerate}
\end{lemma}

\demo   The base orbifold $Q$ has one boundary component and no corners. If $\chi^{\orb}Q>0$ then $Q$ is a disc with at most one cone point thus  $M= D^2\times S^1.$ If $\chi^{\orb}Q=0$ then $Q$  is a Mobius band or a disc with two cone points labeled $2$ and in either case $Q$ has a $2$--fold orbifold-cover that is an annulus $A.$ But then $M$ is $2$--fold covered by a circle bundle over $A.$ Since $M$ is orientable it follows that this bundle is $S^1\times A$ and hence $M$ is a twisted $I$--bundle over the Klein bottle.

Finally, if $\chi^{\orb}(Q)<0$ then all but one filling of $M$ is Seifert fibered.  There are infinitely many fillings of $M$ which give a Seifert fibered space, $P,$ with base orbifold $Q'$ and $\chi^{\orb}(Q')<0.$ There is an orbifold-covering of $Q'$ which is a closed surface of negative Euler characteristic. The induced covering of $P$ contains an essential vertical torus and is therefore  virtually Haken. \endproof
   
\medskip 
 {\bf Definitions}\qua A {\em slope} on a torus $T$ is the isotopy class of  an essential simple closed curve on $T.$ We say that a slope {\em lifts} to a covering of $T$ if it is represented by a loop which lifts. The following is immediate:
\begin{lemma}\label{slopeslift} Suppose $\tilde{T}\rightarrow T$ is a finite covering. Then the following are equivalent:
\begin{enumerate}
\item Some slope on $T$  lifts to $\tilde{T}.$
\item The covering is finite cyclic.
\item Infinitely many slopes on $T$ lift to $\tilde{T}.$
\end{enumerate} 
\end{lemma} 
 
 The {\em distance,} $\Delta(\alpha,\beta),$ between slopes $\alpha,\beta$ on $T$ is the minimum number of intersection points between representative loops. If $\alpha$ is a slope on a torus boundary component of $M$ then $M(\alpha)$ denotes the manifold obtained by Dehn-filling $M$ using $\alpha.$  A surface $S$ in a manifold $M$ is {\em essential}  if it is compact, connected, orientable, incompressible, properly-embedded, and not boundary-parallel.
 Let $M$ be a manifold with  boundary a torus and
 $\alpha\subset\partial M$ a slope. Suppose that $N$ is a finite cover
 of $M.$ An essential surface $S\subset N$ has {\em invariant slope}
 $\alpha$ if $\partial S\ne\phi$ and every component of $\partial S$
 projects to a loop homotopic to a non-zero multiple of $\alpha.$ We
 call a finite  cover $p\co N\rightarrow M$  a {\em $\partial$--cover} if there is an integer $d>0$  and a homomorphism $\theta\co\pi_1(\partial M)\rightarrow{\Bbb Z}_d$ such that for every boundary component $T$ of $N$  we have $p_*(\pi_1T) = \ker\theta.$  The existence of $\theta$ ensures each component of $\partial N$ is the same cyclic cover of $\partial M.$
 
 The following lemma reduces the proof of the main theorem to constructing an essential non-fiber surface of invariant slope in a $\partial$-cover of $M.$ 

 \begin{lemma}\label{mainlemma} Suppose that $M$ is a compact, connected, orientable irreducible 3--manifold with one torus boundary component. Suppose that there is a $\partial$--cover $N$ of $M$ and an essential non-separating surface  $S\subset N$ of invariant slope. Assume  that $S$  is not a
fiber of a fibration of $N.$ Then $M$ has infinitely many
virtually-Haken Dehn-fillings.\end{lemma} 

\demo We first remark that the particular case that concerns us in
this paper is that $M$ is a bundle with boundary and thus $M$ is
irreducible. Since $M$ is irreducible at most 3 fillings give
reducible manifolds (Gordon and Luecke \cite{GordonLueckered}). A
cover of an irreducible manifold is irreducible (Meeks and Yau
\cite{MeeksYau}).  Therefore it suffices to show there are infinitely
many fillings of $M$ which have a finite cover containing an essential
surface.

If $M$ contains an essential torus then this torus remains
incompressible for infinitely many Dehn-fillings by
Culler--Gordon--Luecke--Shalen \cite[Theorem 2.4.2]{CGLS}. If $M$ is
Seifert fibered then by \fullref{SFSfill} either the result holds or
$M=S^1\times D^2$ or is a twisted $I$--bundle over the Klein bottle.
The latter two possibilities do not contain a surface $S$ as in the
hypotheses. By Thurston's hyperbolization theorem we are reduced to
case that $M$ is hyperbolic.

Since $p\co N\rightarrow M$ is a $\partial$--cover  there is $d>0$ such that every component of $\partial N$ is a $d$--fold cover of $\partial M.$  Let $k$ be a positive integer coprime to $d.$ Let $p_k\co\tilde{N}_k\rightarrow N$ be the $k$--fold cyclic cover dual to $S.$ We claim that there is a homomorphism $\theta_k\co\pi_1M\rightarrow {\Bbb Z}_{kd}$  such that every slope in $\ker\theta_k$ lifts to every component of $\partial\tilde{N}_k.$

Assuming this, the filling $M(\gamma)$ of $M$ is covered by a filling,
$\tilde{N}_k(\gamma),$ of $\tilde{N}_k$ if and only if the slope
$\gamma\subset\partial M$ lifts to each component of
$\partial\tilde{N}_k.$ Since $S$ is non-separating, by Wu
\cite[Theorem 5.7]{Wu2}, there is $K>0$ such that if $k\ge K$ then
there is an essential closed surface $F_k\subset\tilde{N}_k$ obtained
by Freedman tubing two lifts of $S.$ We choose such $k$ coprime to
$d.$ By \cite[Theorem 5.3]{Wu2}, there is a finite set of slopes
$\beta_1,\cdots,\beta_n$ on $\partial M$ and $L>0$ so that if
$\gamma\subset\partial M$ is a slope and $\Delta(\gamma,\beta_i) \ge
L$ for all $i$ then the projection of $F_k$ into $M(\gamma)$ is
$\pi_1$--injective. Assuming the claim, there are infinitely many
slopes $\gamma\in \ker\theta_k$ satisfying these inequalities. For
such $\gamma$ the cover $\tilde{N}_k(\gamma)\rightarrow M(\gamma)$
contains the essential surface $F_k.$

It only remains to prove the claim. Let $T$ be a component of $\partial N$ and $\beta\subset T$ be the slope given by $S\cap T.$  Let $\tilde{T}$ be a component of $\partial\tilde{N}_k$ which covers $T.$  The cover $p_k|\co\tilde{T}\rightarrow T $ is cyclic of degree $k'$ some divisor of $k$ (depending only on $|S\cap T|$). Also $\beta$ lifts to this cover. Suppose that a slope $\gamma\subset\partial M$ lifts to a slope $\tilde{\gamma}\subset T.$ It follows that  $\tilde{\gamma}$ lifts to $\tilde{T}$ if $k'$ divides $\Delta(\tilde{\gamma},\beta).$ If this condition is satisfied by some lift, $\tilde{\gamma},$ of $\gamma$ then, since $S$ has invariant slope and $N\rightarrow M$ is a $\partial$--cover,  it is satisfied by every such lift.

Let $\tilde{T}\rightarrow T$ be the $k'$--fold cyclic cover dual to
$\beta.$ Since $k'$ and $d$ are coprime the composite of this cover and
the cyclic $d$--fold cover $T\rightarrow\partial M$ is a cyclic cover
of degree $dk'.$ By \fullref{slopeslift} there are infinitely many
slopes on $\partial M$ which lift to $\tilde{T}.$ Every slope on
$\partial M$ which lifts to $\tilde{T}$ also lifts to every
component of $\partial\tilde{N}_k.$ This proves the claim.\endproof

\proof[Proof of \fullref{maintheorem}] We attempt to construct  $S$ and $N$ as in \fullref{mainlemma}.  The action of the monodromy on $H_1(F;\zt)$ has some finite order $m.$ Therefore there is a finite cyclic $m$--fold cover $W\rightarrow M$ such that $W$ is a bundle with fiber $F$ and the action of the monodromy for $W$ on $H_1(F;{\Bbb Z}_2)$ is trivial. We then have $$H^1(W;\zt) \cong H^1(F;\zt) \oplus H^1(S^1;\zt).$$

Since $F$ has boundary and $F\ne  D^2$ we may choose a non-zero element $\phi = (b,0) \in H^1(F; \zt) \oplus H^1(S^1;
\zt).$ This determines a two-fold cover $\tilde{W}$ of $W.$ Since $F$ has one boundary component, $\phi$ vanishes on $H_1(\partial W;\zt),$ and since $W$ has one boundary component, $\tilde{W}$  has exactly two boundary components 
$T_1$ and $T_2.$ The action of the covering involution, $\tau,$ swaps these tori. In particular $\tilde{W}\rightarrow M$ is a $\partial$--cover.

We claim that there is an essential  surface $S$  in  $\tilde{W}$ such
that $$\tau_*[S] = -[S] \ne 0 \in H_2(\tilde{W},\partial\tilde{W};{\Bbb Z}).$$  
Using real coefficients, all cohomology groups have direct-sum decomposition into $\pm1$ eigenspaces for $\tau^*;$ thus  $H^1(\partial\tilde{W};{\Bbb R}) = V_+\oplus V_-.$ Since $\tau$ swaps $T_1$ and $T_2$ then, with obvious notation, it swaps $\mu_1$ with $\mu_2$ and $\lambda_1$ with $\lambda_2.$
If $\epsilon=\pm 1$ then   $V_{\epsilon}$ has basis
$\{\mu_1+\epsilon\mu_2,\lambda_1+\epsilon\lambda_2\}$ and thus has
dimension $2.$ Let $$K = \Im \left[ \incl^*\co H^1(\tilde{W};{\Bbb R}) \rightarrow H^1(\partial\tilde{W};{\Bbb R})\right].$$
Decompose $K = K_+\oplus K_-.$ We claim that $\dim(K_+)=\dim(K_-)=1.$ Since $\dim(K)=2$ the only other possibilities are that $K_+=V_+$ or $K_-=V_-.$ The intersection pairing on $\partial\tilde{W}$ is dual to the pairing on $H^1(\partial\tilde{W},{\Bbb R})$ given by $<\phi,\psi>\ =\ (\phi\cup\psi)\cap[\partial\tilde{W}].$  This pairing vanishes on $K.$ Since $<\mu_1+\epsilon\mu_2,\lambda_1+\epsilon\lambda_2>\ =\  2<\mu_1,\lambda_1>\ =\pm2,$ 
 the restriction of $<,>$ to each of $V_{\pm}$ is non-degenerate. This contradicts $K=V_{\pm}.$

Choose a  primitive class  $\phi\in H^1(\tilde{W};{\Bbb Z})$ with $\incl^*\phi\in K_-.$ Let $S$  be an essential  oriented surface in $\tilde{W}$ representing the class Poincar\'e dual to $\phi.$ Then $\tau_*[S] = -[S] $ as required.

The $1$--manifold $\alpha_i=T_i\cap\partial S$ with the induced orientation is a $1$--cycle in $\partial \tilde{W}.$ Then $[\partial
S]=[\alpha_1]+[\alpha_2] \in H_1(\partial\tilde{W}).$ Since $T_i$ is a torus all the components of $\alpha_i$ are parallel. Since $\tau(T_1) = T_2$ all components of $\partial S$ project to isotopic loops in $\partial W$ thus $S$ has invariant slope for the cover $\tilde{W}\rightarrow M.$ This gives:

\medskip
\noindent{\bf Case (i)}\qua If  $S$ is not the fiber of a fibration of
$\tilde{W}$ then the result follows from \mbox{\fullref{mainlemma}}.  

\medskip
 Thus we are left with the case that  $S$ is the fiber of a fibration of $\tilde{W}.$
Let $N$ be the $\zt$--universal   covering of $W.$ This is a regular covering and each component of $\partial N$ is a two-fold cover of $\partial W.$  We claim that the composition of coverings $N\rightarrow W\rightarrow M$  is regular. 

Recall that a subgroup $H<G$ is {\em characteristic} if it is preserved by $Aut(G).$ The $\zt$--universal   covering $N\rightarrow W$ corresponds to the characteristic subgroup $\pi_1N < \pi_1W.$ The cover $W\rightarrow M$ is cyclic and so $\pi_1W$ is normal in $\pi_1M.$ A characteristic subgroup of a normal subgroup is normal. Hence $\pi_1N$ is also normal in $\pi_1M.$ This proves the claim.  It follows that $N\rightarrow M$ is a $\partial$--cover.  A pre-image, $\tilde S,$ of  $S$ in $N$ is a fiber 
of a fibration.

\medskip
\noindent{\bf Case (ii)}\qua Suppose the one-dimensional vector space of $H_2(N, \partial N;
\mathbb{R})$ spanned by $[\tilde S]$ is invariant under the group of covering
transformations of $N \rightarrow M.$ 

\medskip
Then, by \fullref{fibering}, $M$ is semi-fibered which contradicts our hypothesis. This completes case (ii). Therefore there is some covering  transformation, $\sigma$,
such that $\sigma_*[\tilde S]\ne\pm[\tilde S].$

Because $\tilde S$ and $\sigma\tilde S$ are fibers, they both meet every boundary component of $N.$ Since $S$ has invariant slope for the cover $N\rightarrow M$ it follows that $\tilde S$ and $\sigma \tilde{S}$ have the same invariant slope for this cover.

\medskip
\noindent{\bf Case (iii)}\qua Suppose $S$ is a fiber and $[\partial \tilde S]\ne\pm\sigma_*[\partial\tilde S]\ 
\in\ H_1(\partial
N).$  

\medskip
Given a boundary component of
$N$, there are integers $a$ and $b$  such that the class $a[\tilde S]
+ b\cdot\sigma_*[\tilde S] \in H_2(N,\partial N)$ is non-zero and represented by 
an essential surface $G$ that
misses this boundary component.   Thus $G$ is not  a fiber of a
fibration.  Clearly $G$ has invariant slope. The result now follows from \fullref{mainlemma} applied to the surface $G$ in the $\partial$--cover $N.$ This completes case (iii). The remaining case is:

\medskip
 \noindent  {\bf Case (iv)}\qua $S$ is a fiber and  there is $\epsilon\in\{\pm1\}$ with $\sigma_*[\partial\tilde 
S]=\epsilon\cdot[\partial\tilde S]\ \in\
H_1(\partial N).$ 

\medskip
Consideration of the homology exact sequence for the pair $(N,\partial N)$ shows $x = \sigma_*[\tilde{S}]  - \epsilon\cdot [\tilde{S} ]  \in H_2(N,\partial N)$ is the image of some $y\in H_2(N).$ Using exactness of the sequence again it follows that $y+i_*H_2(\partial N)$ is not zero in $H_2(N)/i_*H_2(\partial N)$. Hence every filling of $N$ produces a closed manifold with $\beta_2>0.$ Infinitely many slopes on $\partial M$ {\em lift} to slopes on $\partial N.$ The result follows. This completes the proof of case (iv) and thus of the \fullref{maintheorem}. \endproof

\proof[Proof of \fullref{knot}]
Let $\eta(K)$ be an open tubular neighborhood of $k.$ By hypothesis the knot exterior  $M = N\setminus\eta(K)$ is irreducible.  Every semibundle contains two disjoint compact surfaces whose union
 is non-separating, thus the first Betti number with mod-2 coefficients of a semi-bundle is at least $2.$ Because $N$ is a homology sphere  $H_1(M;\zt)\cong\zt,$ therefore $M$ is not a semi-bundle. Since $N$ is a homology sphere it, and therefore $M,$  are orientable.
 
  If $M$ is a bundle with fiber $F$ then, since $N$ is a homology sphere, $F$ has exactly one boundary component. Since $k$ does not bound a disk in $N$ it follows that $M\ne D^2\times S^1.$ The  result now follows from \fullref{maintheorem}. If $M$ contains a closed essential surface then infinitely many fillings are Haken, \cite[Theorem 2.4.2]{CGLS}. The remaining possibilities are that $M$ is hyperbolic and not a bundle, or else Seifert fibered. The hyperbolic non-bundle case follows from \cite{CooperLongvirt}. 
  
This leaves the case that $M$ is Seifert fibered. The manifold $M$ is not a twisted $I$--bundle over the Klein bottle because the latter has mod-2 Betti number $2$. The result now follows from \fullref{SFSfill}.\endproof

\bibliographystyle{gtart}
\bibliography{link}

\begin{thebibliography}{}
\providecommand\bibmarginpar{\leavevmode\marginpar}
\def\urlstyle#1{{\tt #1}}

\bibitem{CooperLongvirt}
\textbf{D Cooper}, \textbf{D\,D Long}, \emph{Virtually {H}aken {D}ehn-filling},
  J. Differential Geom. 52 (1999) 173--187 \xox{MR}{1743462}

\bibitem{CooperWalsh2}
\textbf{D Cooper}, \textbf{G\,S Walsh},
  \href{http://dx.doi.org/10.2140/gt.2006.10.2247} {\emph{Three-manifolds,
  virtual homology, and the group determinant}}, Geom. Topol. 10 (2006)
  2247--2269

\bibitem{CGLS}
\textbf{M Culler}, \textbf{C\,M Gordon}, \textbf{J Luecke}, \textbf{P\,B
  Shalen},
  \href{http://links.jstor.org/sici?sici=0003-486X(198703)2:125:2%3C237:DSOK%3%
E2.0.CO%3B2-7} {\emph{Dehn surgery on knots}}, Ann. of Math. $(2)$ 125 (1987)
  237--300 \xox{MR}{881270}

\bibitem{GordonLueckered}
\textbf{C\,M Gordon}, \textbf{J Luecke},
  \href{http://dx.doi.org/10.1016/0040-9383(95)00016-X} {\emph{Reducible
  manifolds and {D}ehn surgery}}, Topology 35 (1996) 385--409 \xox{MR}{1380506}

\bibitem{Hassmin}
\textbf{J Hass}, \href{http://dx.doi.org/10.1007/BF01161639} {\emph{Surfaces
  minimizing area in their homology class and group actions on 3--manifolds}},
  Math. Z. 199 (1988) 501--509 \xox{MR}{968316}

\bibitem{HempelJaco}
\textbf{J Hempel}, \textbf{W Jaco},
  \href{http://links.jstor.org/sici?sici=0003-486X(197201)2:95:1%3C86:FGO3WA%3%
E2.0.CO%3B2-B} {\emph{Fundamental groups of 3--manifolds which are
  extensions}}, Ann. of Math. $(2)$ 95 (1972) 86--98 \xox{MR}{0287550}

\bibitem{MeeksYau}
\textbf{W Meeks~III}, \textbf{L Simon}, \textbf{S\,T Yau},
  \href{http://links.jstor.org/sici?sici=0003-486X(198211)2:116:3%3C621:EMSES%
A%3E2.0.CO%3B2-R} {\emph{Embedded minimal surfaces, exotic spheres, and
  manifolds with positive {R}icci curvature}}, Ann. of Math. $(2)$ 116 (1982)
  621--659 \xox{MR}{678484}

\bibitem{Wolf}
\textbf{J\,A Wolf}, \emph{Spaces of constant curvature}, third edition, Publish
  or Perish, Boston (1974) \xox{MR}{0343214}

\bibitem{Wu2}
\textbf{Y-Q Wu}, \href{http://dx.doi.org/10.1016/S0040-9383(03)00046-6}
  {\emph{Immersed essential surfaces and {D}ehn surgery}}, Topology 43 (2004)
  319--342 \xox{MR}{2052966}

\bibitem{Z1}
\textbf{L Zulli}, \href{http://dx.doi.org/10.1016/S0166-8641(96)00175-7}
  {\emph{Semibundle decompositions of 3--manifolds and the twisted
  cofundamental group}}, Topology Appl. 79 (1997) 159--172 \xox{MR}{1464194}

\bibitem{Z2}
\textbf{L Zulli}, \emph{Seifert 3--manifolds that are bundles and
  semi-bundles}, Houston J. Math. 27 (2001) 533--540 \xox{MR}{1864797}

\end{thebibliography}

\end{document}